\documentclass[12pt,a4paper]{article}
\usepackage{a4,amsmath,amssymb,amsfonts,amsthm}
\usepackage{amsrefs}
\usepackage[all]{xy}

\begin{document}
\bibliographystyle{plain}
 

\def\mR{\M{R}}           
\def\mZ{\M{Z}}           
\def\mN{\M{N}}           
\def\mQ{\M{Q}}       
\def\mC{\M{C}}  
\def\mG{\M{G}}



\def\Spec{{\rm Spec}}
\def\rg{{\rm rg}}
\def\Hom{{\rm Hom}}
\def\Aut{{\rm Aut}}
 \def\Tr{{\rm Tr}}
 \def\Exp{{\rm Exp}}
 \def\Gal{{\rm Gal}}
 \def\End{{\rm End}}
 \def\det{{{\rm det}}}
 \def\Td{{\rm Td}}
 \def\ch{{\rm ch}}
 \def\che{{\rm ch}_{\rm eq}}
  \def\Spec{{\rm Spec}}
\def\Id{{\rm Id}}
\def\Zar{{\rm Zar}}
\def\Supp{{\rm Supp}}
\def\eq{{\rm eq}}
\def\Ann{{\rm Ann}}
\def\LT{{\rm LT}}
\def\Pic{{\rm Pic}}
\def\rg{{\rm rg}}
\def\et{{\rm et}}
\def\sep{{\rm sep}}
\def\ppcm{{\rm ppcm}}
\def\ord{{\rm ord}}
\def\Gr{{\rm Gr}}
\def\ker{{\rm ker}}
\def\rk{{\rm rk}}


\def\beginProof{\par{\bf Proof. }}
 \def\endProof{${\qed}$\par\smallskip}
 \def\pr{^{\prime}}
 \def\prpr{^{\prime\prime}}
 \def\mtr#1{\overline{#1}}
 \def\ra{\rightarrow}
 \def\mfp{{\mathfrak p}}
 
 \def\mQ{{\Bbb Q}}
 \def\mR{{\Bbb R}}
 \def\mZ{{\Bbb Z}}
 \def\mC{{\Bbb C}}
 \def\mN{{\Bbb N}}
 \def\mF{{\Bbb F}}
 \def\mA{{\Bbb A}}
  \def\mG{{\Bbb G}}
 \def\CI{{\cal I}}
 \def\CE{{\cal E}}
 \def\CJ{{\cal J}}
 \def\CH{{\cal H}}
 \def\CO{{\cal O}}
 \def\CA{{\cal A}}
 \def\CB{{\cal B}}
 \def\CC{{\cal C}}
 \def\CK{{\cal K}}
 \def\CL{{\cal L}}
 \def\CI{{\cal I}}
 \def\CM{{\cal M}}
\def\CP{{\cal P}}
\def\CR{{\cal R}}
\def\CG{{\cal G}}
 \def\wt#1{\widetilde{#1}}
 \def\mod{{\rm mod\ }}
 \def\refeq#1{(\ref{#1})}
 \def\blb{{\big(}}
 \def\brb{{\big)}}
\def\mc{{{\mathfrak c}}}
\def\mcpr{{{\mathfrak c}'}}
\def\mcprpr{{{\mathfrak c}''}}
\def\ss{{\rm ss}}
\def\parf{{\rm parf}}
\def\P1{{{\bf P}^1}}
\def\cod{{\rm cod}}
\def\pr{\prime}
\def\prpr{\prime\prime}
\def\ss{\scriptstyle}
\def\OX{{ {\cal O}_X}}
\def\mpartial{{\mtr{\partial}}}
\def\inv{{\rm inv}}
\def\indlim{\underrightarrow{\lim}}
\def\prolim{\underleftarrow{\lim}}
\def\pprolim{'\prolim'}
\def\Pro{{\rm Pro}}
\def\Ind{{\rm Ind}}
\def\Ens{{\rm Ens}}
\def\without{\backslash}
\def\pbdb{{\Pro_b\ D^-_c}}
\def\qc{{\rm qc}}
\def\Com{{\rm Com}}
\def\an{{\rm an}}
\def\gfield{{\rm\bf k}}
\def\s{{\rm s}}
\def\dR{{\rm dR}}
\def\ari#1{\widehat{#1}}
\def\ul#1{\underline{#1}}
\def\sul#1{\underline{\scriptsize #1}}
\def\mou{{\mathfrak u}}
\def\ich{\mathfrak{ch}}
\def\cl{{\rm cl}}
\def\K{{\rm K}}
\def\R{{\rm R}}
\def\F{{\rm F}}
\def\L{{\rm L}}
\def\pgcd{{\rm pgcd}}
\def\rc{{\rm c}}
\def\N{{\rm N}}
\def\E{{\rm E}}
\def\H{{\rm H}}
\def\CHOW{{\rm CH}}
\def\A{{\rm A}}
\def\d{{\rm d}}
\def\Res{{\rm  Res}}
\def\GL{{\rm GL}}
\def\Alb{{\rm Alb}}
\def\alb{{\rm alb}}
\def\Hdg{{\rm Hdg}}
\def\Num{{\rm Num}}
\def\Irr{{\rm Irr}}
\def\Frac{{\rm Frac}}
\def\Sym{{\rm Sym}}
\def\indlim{\underrightarrow{\lim}}
\def\prolim{\underleftarrow{\lim}}


\def\RHom{{\rm RHom}}
\def\rRHom{{\mathcal RHom}}
\def\rHom{{\mathcal Hom}}
\def\dotimes{{\overline{\otimes}}} 
\def\Ext{{\rm Ext}}
\def\rExt{{\mathcal Ext}}
\def\Tor{{\rm Tor}}
\def\rTor{{\mathcal Tor}}
\def\SP{{\mathfrak S}}

\def\H{{\rm H}}
\def\D{{\rm D}}
\def\Del{{\mathfrak D}}

 \newtheorem{theor}{Theorem}[section]
 \newtheorem{prop}[theor]{Proposition}
 \newtheorem{propdef}[theor]{Proposition-Definition}
 \newtheorem{cor}[theor]{Corollary}
 \newtheorem{lemma}[theor]{Lemma}
 \newtheorem{sublem}[theor]{sub-lemma}
 \newtheorem{defin}[theor]{Definition}
 \newtheorem{conj}[theor]{Conjecture}

 \parindent=0pt
 \parskip=5pt

 \author{Richard PINK\footnote{
 Department of Mathematics, 
ETH Z\"urich, 8092 Z\"urich, SWITZERLAND}\ \ 
   and 
 Damian R\"OSSLER\footnote{D\'epartement de Math\'ematiques, 
B\^atiment 425, 
Facult\'e des Sciences d'Orsay, 
Universit\'e Paris-Sud, 
91405 Orsay Cedex, FRANCE}\\
 {\small with an appendix by B. K\"ock}}
 \title{On the Adams-Riemann-Roch theorem in positive characteristic}
\maketitle
\begin{abstract}
We give a new proof of the Adams-Riemann-Roch theorem for a smooth projective morphism 
$X\to Y$, in the situation where $Y$ is a regular scheme, which is quasi-projective over $\mF_p$. We also partially answer a question of B. K\"ock.
\end{abstract}

 \section{Introduction}
 \label{intro}
 Let  $Y$ be a regular quasi-projective scheme over an affine finite-dimensional 
 noetherian scheme $S$. 
 Let $X$ be a regular scheme and let $f:X\to Y$ be a projective morphism of schemes.
 Let $k\geqslant 2$ be a natural number and $E$ an element of $K_0(X)$. 
 The Adams-Riemann-Roch theorem asserts that 
\begin{equation}
 \psi^k(R^\bullet f_*(E))=R^\bullet f_*\bigl(\theta^k(L_f)^{-1}\otimes\psi^k(E)\bigr)
 \label{arr}
\end{equation}
 in $K_0(Y)[{1\over k}] := K_0(Y) \otimes_\mZ\mZ[{1\over k}]$. 
 The various symbols appearing in this formula are defined as follows.

The Grothendieck group of locally free coherent sheaves on a scheme $Z$ is denoted by $K_0(Z)$, and the Grothendieck group of coherent sheaves on $Z$ is denoted by $K'_0(Z)$.
  The obvious group morphism $K_0(Z)\to K'_0(Z)$ 
 is an isomorphism whenever $Z$ is regular, noetherian and carries an ample invertible sheaf (see 
 \cite[Th. I.9]{Manin-Lectures}). 
 In particular it is an isomorphism in the case $Z=Y$.
 For simplicity we will denote the class in $K_0(Z)$ of a sheaf $E$ again by $E$.
 
  For $f:X\to Y$ as above there is a unique group morphism $K_0(X)\to K'_0(Y)$
  which sends the class of a locally free  coherent sheaf $E$ on $X$ 
  to the class of the linear combination 
 $\sum_{j\geqslant 0}(-1)^j R^j f_*(E)$ of coherent sheaves on $Y$. 
  The composite of this group morphism with 
 the inverse of the isomorphism $K_0(Y) \stackrel{\sim}{\to} K'_0(Y)$ 
 is a group morphism $R^\bullet f _*:K_0(X)\to K_0(Y)$. 
 
 To define the symbol $\psi^k$, recall that the tensor product of $\CO_Z$-modules makes the group $K_0(Z)$  into 
 a commutative unitary ring and that the inverse image of coherent sheaves under any morphism of schemes
 $Z'\to Z$ induces a morphism of unitary rings 
 $K_0(Z)\to K_0(Z')$ (see \cite[Par. 1]{Manin-Lectures}). Thus $K_0(\cdot)$ may be viewed as 
 a contravariant functor from the category of schemes to the category of 
 commutative unitary rings. 
 The symbol $\psi^k$ refers to an endomorphism 
 of this functor (sic!) that is uniquely determined by the further property that 
$$
 \psi^k(L)=L^{\otimes k}
 $$
for any invertible sheaf $L$ (see \cite[Par. 16]{Manin-Lectures}).
 
The symbol $\theta^k$ refers to a different operation associating 
an element of  $K_0(Z)$ to any locally free coherent sheaf on $Z$. 
It is uniquely determined by the properties:
 \begin{description}
 \item[(i)] for any invertible sheaf $L$ on $Z$ we have
 $$
 \theta^k(L)=1+L+\dots+L^{k-1} ,
 $$
 \item[(ii)] for any short exact sequence
 $
 0\to E'\to E\to E''\to 0
 $
 of locally free coherent sheaves on $Z$ we have
 $$
 \theta^k(E)=\theta^k(E')\otimes\theta^k(E'') ,
 $$
 \item[(iii)] 
 for any morphism of schemes $g:Z'\to Z$ 
and any locally free coherent sheaf $E$ on $Z$
  we have
 $$
 g^*(\theta^k(E))=\theta^k(g^*(E)).
 $$
 \end{description}
 If $Z$ is quasi-projective over an affine finite-dimensional noetherian scheme, 
 it is known that 
 $\theta^k(E)$ is invertible in $K_0(Z)[{1\over k}]$ 
for every locally free coherent sheaf $E$ on $Z$ 
 (see \cite[Lemma 4.3]{Koeck-Grothendieck}). 
  In that case $\theta^k$ extends to a unique map
 $K_0(Z) \to K_0(Z)[{1\over k}]$ satisfying
 $$\theta^k(E)=\theta^k(E')\cdot \theta^k(E'')$$
 whenever $E = E'+E''$ in $K_0(Z)$.

 The symbol $L_f$ denotes the relative cotangent complex of the morphism 
 $f\!:X\to\nobreak Y$. Since $f$ is by construction a local complete intersection morphism, 
 its cotangent complex can be represented by a bounded complex of coherent sheaves 
 (see for instance \cite[Illusie, chap. 8, App. G., 8.5.29]{Fantechi-FGA}). Thus 
 $L_f$ determines a unique element of $K'_0(X)=K_0(X)$, for which 
$\theta^k(L_f) \in K_0(X)[{1\over k}]$ is well-defined by additivity. 
For example, if $f$ is smooth, then $\theta^k(L_f) = \theta^k(\Omega_f)$. 

This explains all the ingredients of the formula \refeq{arr}.
  
The formula \refeq{arr} is classically proven using deformation to 
the normal cone and considering closed immersions and 
relative projective spaces separately (see \cite{BFM}). 

Our aim in this text is to provide a new and more direct proof of 
the formula \refeq{arr} in the specific situation where 
$k$ is a prime number $p$, $f$ is smooth 
and $S$ is a scheme of characteristic $p$, which is of finite type over the finite field $\mF_p$.
 
The search for this proof was motivated by the fact that for any scheme 
$Z$ of characteristic $p$,  the endomorphism $\psi^p:K_0(Z)\to K_0(Z)$ 
coincides with 
the endomorphism $F^*_Z:K_0(Z)\to K_0(Z)$ induced by pullback by the 
absolute Frobenius endomorphism $F_Z:Z\to Z$. This is a consequence of 
 the splitting principle \cite[Par. 5]{Manin-Lectures}.  
We asked ourselves whether in this case $\theta^p(L_f)=\theta^p(\Omega_f)$ 
can also be represented by an explicit virtual bundle.
If such a representative were available, one might try to give a 
direct proof of \refeq{arr} that does not involve factorisation. 
The proof given in Section \ref{arrproof} shows that this is indeed possible. 

In the article 
\cite[sec. 5]{Koeck-Riemann} by B. K\"ock, a different line of speculation lead to 
a question (Question 5.2) in the context of a characteristic $p$ interpretation of the 
Adams-Riemann-Roch formula. Our Proposition \ref{propE} and Proposition \ref{lemmaF} show 
that the answer to this question is positive for a large class of morphisms. See the end of section \ref{arrproof} for details. 

Fix $k\geqslant 2$ and suppose that $Y$ is the spectrum of a finite field. 
The formula \refeq{arr} then formally 
implies the Hirzebruch-Riemann-Roch theorem for $X$ over that field. 
This 
is explained for instance in \cite[Intro.]{Nori-The-Hirzebruch}. 
On the other hand, a specialization argument 
shows that the Hirzebruch-Riemann-Roch theorem for 
varieties over any field follows from the Hirzebruch-Riemann-Roch theorem 
for varieties over finite fields. Thus by reduction modulo primes our proof 
of \refeq{arr} in positive characteristic leads to a proof of 
the Hirzebruch-Riemann-Roch formula in general.

The structure of the article is the following. In Section \ref{bundrep}, 
we construct a canonical bundle representative for the 
element $\theta^p(E)$ for any locally free coherent sheaf $E$ on 
a scheme of characteristic $p$. In Section \ref{arrproof}, we give 
the computation proving \refeq{arr} in the situation where $k=p$, $f$ is smooth 
and $S$ is a scheme of characteristic $p$, which is of finite type over 
$\mF_p$.

\section{A bundle representative for $\theta^p(E)$}
\label{bundrep}

Let $p$ be a prime number and $Z$ a scheme of characteristic $p$.
Let $E$ be a locally free coherent sheaf on $Z$. For any integer 
$k\geqslant0$ let $\Sym^k(E)$ denote the $k$-th symmetric power of $E$.
Then 
$$\Sym(E) := \bigoplus_{k\geqslant0} \Sym^k(E)$$
is a quasi-coherent graded $\CO_Z$-algebra, called the symmetric algebra of $E$.
Let $\CJ_E$ denote the graded sheaf of ideals of $\Sym(E)$ that is locally 
generated by the sections $e^p$ of $\Sym^p(E)$ for all sections $e$ of~$E$, and set
$$\tau(E) := \Sym(E) / \CJ_E.$$
Locally this construction means the following. 
Consider an open subset $U\subset Z$ such that $E|U$ is free, and 
choose a basis $e_1,\dots, e_r$. Then $\Sym(E)|U$ is the polynomial 
algebra over $\CO_Z$ in the variables $e_1,\dots, e_r$. Since $Z$ has
characteristic~$p$, for any open subset $V\subset U$ and any sections 
$a_1,\dots, a_r \in \CO_Z(V)$ we have
$$\bigl(a_1e_1+\ldots+a_re_r\bigr)^p = a_1^pe_1^p+\ldots+a_r^pe_r^p.$$
It follows that $\CJ_E|U$ is the sheaf of ideals of $\Sym(E)|U$ 
that is generated by $e_1^p,\dots, e_r^p$. Clearly that description
is independent of the choice of basis and compatible with localization; 
hence it can be used to an equivalent definition of $\CJ_E$ and $\tau(E)$.

The local description also implies that $\tau(E)|U$ is free over $\CO_Z|U$ 
with basis the images of the monomials $e_1^{i_1}\cdots e_r^{i_r}$ for all 
choices of exponents $0\leqslant i_j<p$. From this we deduce:

\begin{lemma}
\label{lemmaG}
If $E$ is a locally free coherent sheaf of rank $r$, 
then $\tau(E)$ is a locally free coherent sheaf of rank $p^r$. 
\end{lemma}

Now we go through the different properties that characterize the operation~$\theta^p$. 

\begin{lemma}
\label{lemtau(i)}
For any invertible sheaf $L$ on $Z$ we have
 $$\tau(L) \cong \CO_Z \oplus L \oplus\dots\oplus L^{\otimes(p-1)}.$$
\end{lemma}

\beginProof
In this case the local description shows that $\CJ_L$ is the sheaf of ideals 
of $\Sym(L)$ that is generated by $\Sym^p(L)= L^{\otimes p}$. 
The lemma follows at once.
\endProof

\begin{lemma}
\label{lemtau(iii)}
For any morphism of schemes $g:Z'\to Z$ and any locally free coherent sheaf $E$ on $Z$
we have
$$g^*(\tau(E)) \cong \tau(g^*(E)).$$
\end{lemma}

\beginProof
Direct consequence of the construction.
\endProof

\begin{lemma}
\label{lemtau(pre-ii)}
For any two locally free coherent sheaves $E'$ and $E''$ on $Z$ we have
 $$\tau(E'\oplus E'') \cong \tau(E')\otimes\tau(E'').$$
\end{lemma}

\beginProof
The homomorphism of sheaves
$$E'\oplus E'' \hookrightarrow \Sym(E')\otimes\Sym(E''),\ 
  (e',e'') \mapsto e'\otimes 1 + 1\otimes e''$$
induces an algebra isomorphism
$$\Sym(E'\oplus E'') \to \Sym(E') \otimes \Sym(E'').$$
The local description as polynomial rings in terms of bases of $E'|U$ and 
$E''|U$ shows that this is an isomorphism of sheaves of $\CO_Z$-algebras.
Since 
$$(e'\otimes 1 + 1\otimes e'')^p = 
  e^{\prime p}\otimes 1 + 1\otimes e^{\prime\prime p}$$
for any local sections $e'$ of $E'$ and $e''$ of $E''$,
this isomorphism induces an isomorphism of sheaves of ideals 
$$\CJ_{E'\oplus E''} \to
  \CJ_{E'}\otimes\Sym(E'') \oplus \Sym(E')\otimes\CJ_{E''}.$$
The lemma follows from this by taking quotients.
\endProof

\begin{lemma}
\label{lemtau(ii)}
For any short exact sequence $0\to E'\to E\to E''\to 0$ of locally free coherent 
sheaves on a noetherian scheme $Z$ we have
$$\tau(E) = \tau(E')\otimes\tau(E'')$$
in $K_0(Z)$.
\end{lemma}

\beginProof
Let $\wt{E}'$ and $\wt{E}''$ denote the inverse images of $E'$ and $E''$ 
under the projection morphism $Z\times\P1\to Z$. Then there exists 
a short exact sequence
$$0\to \wt{E}'\to\wt{E}\to\wt{E}''\to 0$$
of locally free coherent sheaves on $Z\times\P1$ whose restriction to 
the fiber above $0\in\P1$ is the given one and whose restriction to the 
fiber above $\infty\in\P1$ is split (the construction is given in \cite[I, Par. f)]{Bismut-Gillet-Soule-Analytic}). Thus the respective restrictions 
satisfy $\wt{E}_0\cong E$ and $\wt{E}_\infty\cong E'\oplus E''$.
Using Lemmata \ref{lemtau(iii)} and \ref{lemtau(pre-ii)} this implies that
$$\tau(E) \cong \tau(\wt{E}_0) \cong \tau(\wt{E})_0$$
and
$$\tau(E')\otimes\tau(E'') \cong \tau(E'\oplus E'') 
\cong \tau(\wt{E}_\infty)\cong \tau(\wt{E})_\infty.$$
But the fact that $K_0(Z\times\P1)$ is generated by 
the powers of ${\cal O}(1)$ over $K_0(Z)$ (see \cite[Par. 5]{Manin-Lectures}) 
implies that the restriction to $0$ and $\infty$ induce the same map
$K_0(Z\times\P1) \to K_0(Z)$. Thus it follows that 
$\tau(\wt{E})_0 = \tau(\wt{E})_\infty$ in $K_0(Z)$, whence the lemma.
\endProof

{\bf Remark.} Lemma \ref{lemtau(ii)} can also be proved by an explicit 
calculation of sheaves. For a sketch consider the decreasing filtration of 
$\Sym(E)$ by the graded ideals $\Sym^i(E') \cdot \Sym(E)$ for all $i\geqslant0$. 
One first shows that the associated bi-graded algebra is isomorphic
to $\Sym(E')\otimes\Sym(E'')$. The filtration of $\Sym(E)$ also induces
a filtration of $\tau(E)$ by graded ideals, whose associated bi-graded 
algebra is therefore a quotient to $\Sym(E')\otimes\Sym(E'')$. To prove
that this quotient is isomorphic to $\tau(E')\otimes\tau(E'')$ one shows
that the kernel of the quotient morphism
$\Sym(E')\otimes\Sym(E'') \twoheadrightarrow {\rm Gr}(\tau(E))$
is precisely $\CJ_{E'}\otimes\Sym(E'') \oplus \Sym(E')\otimes\CJ_{E''}$.
But this is a purely local assertion, for which one can assume that the 
exact sequence splits. The calculation then becomes straightforward, 
as in Lemma \ref{lemtau(pre-ii)}.

\begin{prop}
\label{propE}
For any locally free coherent sheaf $E$ on a noetherian scheme $Z$ we have
$\tau(E) = \theta^p(E)$ in $K_0(Z)$.
\end{prop}

\beginProof
Combination of Lemmata \ref{lemtau(i)}, \ref{lemtau(iii)}, \ref{lemtau(ii)}
and the defining properties (i), (ii), (iii) of $\theta^p(\cdot)$
in Section \ref{intro}.
\endProof

\section{Proof of the Adams-Riemann-Roch formula}
\label{arrproof}

Let us now consider the morphism $f:X\to Y$ of the introduction. 
Recall that $Y$ is regular and quasi-projective over an affine
noetherian finite-dimensional scheme $S$ and that $f$ is projective. We make the 
supplementary hypothesis 
that $f$ is smooth and that $S$ is a scheme of characteristic $p$, which 
is of finite type over the finite field $\mF_p$. To prove the formula \refeq{arr} we may also 
suppose that $Y$ and $X$ are connected and thus integral. Then 
$f$ has constant fibre dimension, say $r$.
 
Consider the commutative diagram
$$
\xymatrix{
X \ar[r]^{F}\ar[dr]^{f}\ar@/^2pc/[rr]^{F_X} & X'\ar[r]^{J}\ar[d]^{f'} & X\ar[d]^f\\
                            & Y\ar[r]^{F_Y} & Y
}
$$
where $F_X$ and $F_Y$ are the respective absolute Frobenius morphisms 
and the square is cartesian. The morphism $F=F_{X/X'}$ is called the 
relative Frobenius morphism of $X$ over~$Y$. The following lemma 
summarizes the properties of $F$ that we shall need. For its proof, see
 \cite[Th. 15.7]{Kunz-Kaehler}. 
\begin{lemma}
\label{lemmaFrob}
The morphism $F$ is finite and flat of constant degree $p^r$.
\end{lemma}

Let $I$ denote the kernel of the natural morphism of $\CO_X$-algebras 
$F^*F_*\CO_X\to \CO_X$, which by construction is a sheaf of ideals of
$F^*F_*\CO_X$. Let 
$$\Gr(F^*F_*\CO_X) := \bigoplus_{k\geqslant 0}I^k/I^{k+1}$$
denote the associated graded sheaf of $\CO_X$-algebras.
Let $\Omega_f$ denote the relative sheaf of differentials of~$f$.

\begin{prop}
\label{lemmaF}
There is a natural isomorphism of $\CO_X$-modules
$$I/I^2\cong\Omega_f$$
and a natural isomorphism of graded $\CO_X$-algebras
$$\tau(I/I^2)\cong\Gr(F^*F_*\CO_X).$$
\end{prop}

\beginProof
Since $F$ is affine (see Lemma \ref{lemmaFrob}), there is a canonical 
isomorphism
$$\Spec\ F^*F_*\CO_X \cong X\times_{X'}X,$$
for which the natural morphism of $\CO_X$-algebras $F^*F_*\CO_X\to \CO_X$
corresponds to the diagonal embedding $X\hookrightarrow X\times_{X'}X$.
We carry out these identifications throughout the remainder of this proof.
Then $I$ is the sheaf of ideals of the diagonal, and so $I/I^2$ is 
naturally isomorphic to the relative sheaf of differentials $\Omega_F$.
On the other hand we have $F^*\Omega_{f'} = F^*J^* \Omega_f = 
F_X^* \Omega_f$, which yields a natural exact sequence 
$$F_X^* \Omega_f 
\to\Omega_f\to\Omega_F\to 0.$$
Here the leftmost arrow sends any differential $dx$ to $d(x^p) = 
p\cdot x^{p-1}\cdot dx = 0$. Thus the exact sequence yields an isomorphism 
$\Omega_f\cong\Omega_F \cong I/I^2$, proving the first assertion.

For the second assertion observe that, by the universal property of the 
sym\-metric algebra $\Sym(\cdot)$, the embedding $I/I^2 \hookrightarrow 
\Gr(F^*F_*\CO_X)$ extends to a unique morphism of $\CO_X$-algebras
$$\rho:\Sym(I/I^2)\to\Gr(F^*F_*\CO_X).$$
We want to compare the kernel of $\rho$ with $\CJ_{I/I^2}$. For this recall
that $I$, as the sheaf of ideals of the diagonal, is generated by the 
sections $s\otimes1-1\otimes s$ for all local sections $s$ of $\CO_X$. 
The $p$-th power of any such section is 
$$(s\otimes1-1\otimes s)^p = s^p\otimes1-1\otimes s^p = 0$$
in $F^*F_*\CO_X$, because $s^p = F_X^*s$ is the pullback via $F_X$ of a section 
of $\CO_X$ and hence also the pullback via $F$ of a section of $\CO_{X'}$. 
Thus $\rho$ sends the $p$-th powers of certain local generators of $I/I^2$ to zero. 
But in Section \ref{bundrep} we have seen that $\CJ_{I/I^2}$ is locally
generated by the $p$-th powers of any local generators of $I/I^2$. Therefore 
$\rho(\CJ_{I/I^2})=0$, and so $\rho$ factors through a morphism of $\CO_X$-algebras
$$\bar\rho: \tau(I/I^2)\to\Gr(F^*F_*\CO_X).$$
From the definition of $\Gr(F^*F_*\CO_X)$ we see that $\rho$ and 
hence $\bar\rho$ is surjective.

On the other hand the smoothness assumption on $f$ implies that 
$I/I^2\cong\Omega_f$ is locally free of rank~$r$. Thus Lemma 
\ref{lemmaG} shows that $\tau(I/I^2)$ is locally free of rank~$p^r$. 
By Lemma \ref{lemmaFrob} the same is true for $F^*F_*\CO_X$ and hence,
since $X$ is integral, for $\Gr(F^*F_*\CO_X)$ at the generic point $\eta$
of~$X$. As $\bar\rho$ is surjective, it is therefore an isomorphism
at $\eta$. Therefore the sheaf $\ker(\bar\rho)$ vanishes at $\eta$.
But since $X$ is integral, any torsion subsheaf of a locally free coherent 
sheaf on $X$ is zero. Thus $\ker(\bar\rho)=0$ everywhere, and so $\bar\rho$
is the desired isomorphism.
\endProof

{\bf Remark.} The assumption that $f$ is projective was not used in the proof of 
Proposition \ref{lemmaF}. In particular, its conclusion is valid 
without this assumption.

\begin{lemma}
\label{lemmaH}
Let $Z$ be a quasi-projective scheme of finite dimension over an affine 
noetherian scheme.
Let $E$ be a locally free coherent sheaf of rank $r$ on $Z$.
Then the class of $E$ is invertible in the ring $K_0(Z)[{1\over r}]$.
\end{lemma}
\beginProof
Let $F^1 K_0(Z)$ be the kernel of the rank morphism $K_0(Z)\to\mZ$. 
This is an ideal whose $k$-th power vanishes for all $k>\dim(Z)$
(for this see \cite[V, par. 3, Cor. 3.10]{Fulton-Lang-Riemann-Roch}). 
The infinite sum in $K_0(Z)[{1\over r}]$ 
$$1/r+(r-E)/r^2+(r-E)^{\otimes 2}/r^3+\dots$$
therefore only has a finite number of non-vanishing terms. 
A direct calculation with geometric series shows that this 
sum is an inverse of $E$ in $K_0(Z)[{1\over r}]$.
\endProof

{\bf Remark.}  In \cite[Question 5.2]{Koeck-Riemann}, B. K\"ock in particular asks the following 
question: is the equation
$$
F_*(\theta^p(\Omega_f)^{-1})=1
$$
valid in $K_0(Y)[{1\over p}]$ ? 
Proposition \ref{lemmaF} implies that the answer to this question is 
positive. Indeed, using the projection formula in $K_0$-theory, we compute
$$
F_*(\theta^p(\Omega_g)^{-1})=F_*((F^*F_*\CO_Z)^{-1})=
F_*(F^*(F_*\CO_Z)^{-1})=(F_*\CO_Z)\otimes(F_*\CO_Z)^{-1}=1.
$$
This computation is partially repeated below.

We now come to the proof of the Adams-Riemann-Roch formula, which 
results from the following calculation in $K_0(X)[{1\over p}]$. This calculation 
is in essence already in \cite[Prop. 5.5]{Koeck-Riemann}. It did not lead to 
a proof of the formula \refeq{arr} there, because the Proposition \ref{lemmaF} was 
missing. 
\begin{eqnarray*}
    \psi^p(R^\bullet f_*(E)) 
&=& F_Y^*R^\bullet f_*(E) \\
&=& R^\bullet f'_*(J^*(E)) \\
&=& R^\bullet f'_*\bigl((F_*\CO_X)\otimes(F_*\CO_X)^{-1}\otimes J^*(E)\bigr) \\
&=& R^\bullet f'_*F_*\bigl(F^*(F_*\CO_X)^{-1}\otimes F^*J^*(E)\bigr) \\
&=& R^\bullet f_*\bigl((F^*F_*\CO_X)^{-1}\otimes F_X^*(E)\bigr)\\
&=& R^\bullet f_*\bigl(\theta^p(\Omega_f)^{-1}\otimes\psi^p(E)\bigr).
\end{eqnarray*}
Here the first equality uses the fact that $\psi^p=F_Y^*$ in $K_0(Y)$.
The second equality follows from the fact that the formation of cohomology 
commutes with flat base change. The third equality is the definition of
$(F_*\CO_X)^{-1}$ in $K_0(X')[{1\over p}]$ using Lemmata \ref{lemmaFrob}
and \ref{lemmaH}. The fourth equality is justified by the projection formula 
in $K_0$-theory (see \cite[Prop. 7.13]{Manin-Lectures}). The fifth equality is just a simplification.
Finally, Proposition \ref{lemmaF} and Proposition \ref{propE} imply that 
$$F^*F_*\CO_X = \Gr(F^*F_*\CO_X) = \tau(I/I^2) =
  \theta^p(I/I^2) = \theta^p(\Omega_f) = \theta^p(L_f)$$
as elements of $K_0(X)$. This and the fact that $\psi^p=F_X^*$ in $K_0(X)$
prove the last equality, and we are done.

\appendix

\begin{flushleft}
\bf Appendix : another formula for the Bott element
\end{flushleft}
\begin{center}
by Bernhard K\"ock\footnote{School of Mathematics, University of Southampton,
SO17 1BJ, United Kingdom. e-mail: B.Koeck@soton.ac.uk}
\end{center}

The object of this appendix is to give a new formula for the Bott
element of a smooth morphism. This formula is analogous to the
final displayed formula in the main part of this paper and
completes a list of miraculous analogies explained in section~5
of~\cite{Koeck-Riemann}; it also streamlines the proof of Theorem~3.1 in section~3
of~\cite{Koeck-Riemann}. It is probably needless to say that this appendix is
inspired by the elegant approach to the Adams-Riemann-Roch theorem
in positive characteristic developed by Richard Pink and Damian
R\"ossler in the main part of this paper.

We begin by setting up the context. Let $l$ be a prime and let
$f:X \ra Y$ be a smooth quasi-projective morphism between
Noetherian schemes of relative dimension $d$. We furthermore
assume that there exists an ample invertible $\CO_X$-module. Let
$\Omega_f$ denote the lcoally free sheaf of relative differentials
and let $\theta^l(\Omega_f) \in K_0(X)$ denote the $l$-th Bott
element associated with $\Omega_f$ (see Introduction). Furthermore
let $\Delta: X \ra X^l$ denote the diagonal morphism from $X$ into
the $l$-fold cartesian product $X^l := X\times_Y \ldots \times_Y
X$. We view $\Delta$ as a $C_l$-equivariant morphism where the
cyclic group~$C_l$ of order $l$ acts trivially on $X$ and by
permuting the factors on $X^l$. In particular we have a pull-back
homomorphism $\Delta^*: K_0(C_l, X^l) \ra K_0(C_l, X)$ between the
corresponding Grothendieck groups of equivariant locally free
sheaves on $X^l$ and $X$, respectively. As the closed immersion
$\Delta$ is also regular we furthermore have a push-forward
homomorphism $\Delta_*: K_0(C_l, X) \ra K_0(C_l, X^l)$ (see
section~3 in~\cite{Koeck-Grothendieck}). Let finally $([\CO_X[C_l]])$ denote the
principal ideal of $K_0(C_l,X)$ generated by the regular
representation $[\CO_X[C_l]]$. We have a natural map $K_0(X) \ra
K_0(C_l,X) \ra K_0(C_l,X)/([\CO_X[C_l]])$ which is in fact
injective under certain rather general assumption (see
Corollary~4.4 in \cite{Koeck-Riemann}). The following theorem strengthens
Theorem~3.1 in~\cite{Koeck-Riemann}; it should be viewed as an analogue of the
formula $\theta^p(\Omega_f) = F^* F_*(\CO_X)$ proved at the very
end of the main part of this paper.

{\bf Theorem.} We have
\[\theta^l(\Omega_f) = \Delta^* (\Delta_* (\CO_X))) \quad \textrm{in }
K_0(C_l, X)/([\CO_X[C_l]]).\]

{\em Proof.} Let $\CI_\Delta$ denote the ideal sheaf corresponding
to the regular closed immersion $\Delta: X \ra X^l$. Then we have
\[\Delta^*(\Delta_*(\CO_X)) =
\lambda_{-1}(\CI_\Delta/\CI_\Delta^2) \quad \textrm{in }
K_0(C_l,X)\] by the equivariant self-intersection formula (see
Corollary~(3.9) in~\cite{Koeck-Grothendieck}); here $\lambda_{-1}(\CE)$ denotes the
alternating sum $[\CO_X] - [\CE] + [\Lambda^2(\CE)] \pm \ldots$
for any locally free $C_l$-sheaf $\CE$ on $X$. Furthermore we know
that $\CI_\Delta/\CI_\Delta^2$ is $C_l$-isomorphic to $\Omega_f
\otimes \CH_{X,l}$ where $\CH_{X,l} := \textrm{ker}(\CO_X[C_l]
\,\, \stackrel{{\rm sum}}{\longrightarrow} \CO_X)$ denotes the
augmentation representation (see Lemma~3.5 in~\cite{Koeck-Riemann}). Finally we
have $\lambda_{-1}(\CE \otimes \CH_{X,l}) = \theta^l(\CE)$ in
$K_0(C_l,X)/([\CO_X[C_l]])$ for any locally free $C_l$-module
$\CE$ on $X$ (see Proposition~3.2 and Remark~3.9 in~\cite{Koeck-Riemann}). Putting
these three facts together we obtain the desired equality of
classes in $K_0(C_l, X)/([\CO_X[C_l]])$.

{\em Remarks.}\\
(a) As in the remark after Lemma~3.3 in the main part of this
paper, using the projection formula, we can easily derive the
original formula $\Delta_*\left(\lambda_{-1}(\Omega_f \otimes
\CH_{X,l})^{-1}\right) =1$ in $K_0(C_l,X)[l^{-1}]/([\CO_X[C_l]])$
(see Theorem~3.1 in \cite{Koeck-Riemann}) from (the proof of) the above theorem. \\
(b) The following table summarizes the astounding analogies
mentioned at the beginning of this appendix. While the left hand
column refers to the situation of the main part of this paper the
right hand column refers to the situation of this appendix and of
section~4 in \cite{Koeck-Riemann}. The entries in the table are of a very symbolic
nature; more detailed explanations can be found in section~5 of
\cite{Koeck-Riemann}. For instance, $\tau^l: K_0(X) \ra K_0(C_l, X)$ and
$\tau^l_{\rm ext}: K_0(X) \ra K_0(C_l, X^l)$ denote the $l$-th
tensor-power operation and $l$-th external-tensor-power operation,
respectively.

\begin{center}
\begin{tabular}{|c||c|} \hline
$\psi^p = F_X^*$ & $\psi^l = \tau^l$\\ \hline relative Frobenius
$F:X \ra X'$ & diagonal $\Delta: X \ra X^l$ \\ \hline
$f$ is smooth & $f$ is smooth \\
$\Rightarrow$ $F$ is flat & $\Rightarrow$ $\Delta$ is regular \\
$\Rightarrow$ We have $F_*:K_0(X) \rightarrow K_0(X')$ &
$\Rightarrow$ We have $\Delta_*:K_0(C_l,X) \rightarrow K_0(C_l,
X^l)$\\ \hline $f': X' \ra Y$ & $f^l: X^l \ra Y$ \\ \hline $J^*:
K_0(X) \ra K_0(X')$ & $\tau_{\rm ext}^l: K_0(X) \ra K_0(C_l,
X^l)$\\ \hline Base change: $F_Y^* f_* =
(f')_*J^*$ & K\"unneth formula: $\tau^l f_* = f^l_* \tau_{\rm ext}^l$\\
\hline
$F_X^* = F^* J^*$ & $\tau^l = \Delta^* \tau^l_{\rm ext}$\\
\hline $\theta^p(\Omega_f) = F^*(F_*(\CO_X))$ &
$\theta^l(\Omega_f) = \Delta^* (\Delta_*(\CO_X))$
\\ \hline
\end{tabular}
\end{center}

\end{document}